\documentclass{amsart}

\usepackage{amsfonts,amsmath,amssymb}

\DeclareMathOperator{\Span}{Span}
\DeclareMathOperator{\Gr}{Gr}
\DeclareMathOperator{\Lie}{Lie}

\DeclareMathOperator{\Iso}{Iso}
\newtheorem{definition}{Definition}
\newtheorem{theorem}{Theorem}
\newtheorem{proposition}{Proposition}
\newtheorem{remark}{Remark}

\newcommand{\is}{_{\text{\rm iso}}}
\newcommand{\iso}{_{\text{\emph{iso}}}}

 \input xy
\xyoption {all}

\newcommand{\p}{^{\prime}}

\newcommand{\A}{\mathcal{A}}
\newcommand{\B}{\mathcal{B}}
\newcommand{\W}{\mathcal{W}}
\newcommand{\T}{\mathcal{T}}

\newcommand{\M}{\mathcal{M}}

\begin{document}
\markboth{R. Lipyanski, N. Vanetik}{On Borel complexity of the isomorphism problems}

\title{ON BOREL COMPLEXITY OF THE ISOMORPHISM PROBLEMS FOR GRAPH RELATED CLASSES OF LIE ALGEBRAS AND FINITE $p$-GROUPS}


\author{\footnotesize RUVIM LIPYANSKI}
\address{Department of Mathematics, Ben Gurion University\\ Beer Sheva, Israel}
\email{lipyansk@math.bgu.ac.il}

\author{\footnotesize NATALIA VANETIK}
\address{Shamoon College of Engineering\\ Beer Sheva, Israel}
\email{natalyav@sce.ac.il}


\begin{abstract}
We reduce the isomorphism problem for undirected graphs without loops  to the isomorphism problems for a class of finite dimensional $2$-step nilpotent Lie algebras over a field, and for a class of finite $p$-groups. We show that the isomorphism problem for graphs is harder than the two latter isomorphism problems in the sense of Borel reducibility. A computable analogue of Borel reducibility was introduced by S. Coskey, J.D. Hamkins, and R. Miller in \cite{Cos}. A relation of the isomorphism problem for undirected graphs to the well-known problem of classifying pairs of matrices over a field (up to similarity) is also studied.

\end{abstract}

\maketitle


\section{Introduction}
 In this paper we define the class $\mathsf{GLA}$ of graph Lie algebras over a field and the class $\mathsf{GpG}$ of graph $p$-groups (Sections 2 and 3, respectively). We also reduce the isomorphism problems for the class $\mathsf{GRAPH}$ of undirected graphs without loops to  the isomorphism problems for the above classes.

Previously, the graph isomorphism problem was reduced to the isomorphism problems for rings \cite{Sax}, algebras \cite{Kim} and groups {\cite {PBel,D}.  Contrary to the paper \cite{Kim}, where algebras were infinite dimensional, we reduce the graph isomorphism problem to the isomorphism problem of a class of 2-step nilpotent finite dimensional Lie algebras over a field. Also, in \cite{D} the isomorphism problem for graphs was reduced to the isomorphism problem for a class of infinite groups.

  A reduction of the graph isomorphism problem to the isomorphism problem for the class $\mathsf{GpG}$ of finite $p$-groups was  given in \cite{PBel}. In Section 3, we present a new proof of this result based on the Lazard correspondence between the category of nilpotent Lie rings of nilpotency class $c$ and order $p^n$, $p>c$, and the category of finite $p$-group of order $p^n$ and nilpotency class $c$.

To compare the complexity of the isomorphism problems for the classes $\mathsf{GLA}$, $\mathsf{GpG}$, and $\mathsf{GRAPH}$ we use the \emph{polynomial time Borel reducibility of equivalence relations} on countable sets (see \cite{Cos}).

Let $A$ and $C$ be two countable sets and $R,S$ be equivalence relations on $A$ and $C$, respectively. We say that $(A,R)$ is computably Borel-reducible to $(C,S)$, and write $(A,R)\leq_{B}^P(C,S)$, if there exists a polynomially computable map $f:A\rightarrow C$, such that for all  $x$ and $y$ in $A$
    $$
    xRy\Leftrightarrow f(x) S f(y).
    $$
  In other words, the reduction function $f$ yields a classification of the elements  of $A$ up to $R$ using invariants from $C/S$.
   We will also say that the classification problem of the elements  of $A$ up to $R$ is  \emph{not harder} (in the Borel sense) than the classification problem of the elements  of $C$ up to $S$. We say that $(A,R)$  and $(C,S)$ are \emph{Borel equivalent} and write $(A, R)\equiv_B^P (C, S)$ if they are polynomial-time Borel-reducible ($P$-Borel-reducible) one to another, i.e., $(A,T)\leq_B^P (C,S)$ and $(C,S)\leq_B^P(A,R)$. A detailed discussion of Borel reducibility is given in Section 4.

  Let $\mathsf{D}$  and $\mathsf{D^\prime}$ be two classes of finite structures and $\mbox{Iso}_\mathsf{D}$, with
  $\mbox{ Iso}_\mathsf{D^\prime}$ beiing two isomorphism relations on these classes, respectively. Then $P$-Borel reducibility of the pair $(\mathsf{D},\rm{Iso}_\mathsf{D})$ to $(\mathsf{D^\prime}, \rm{Iso}_\mathsf{D^\prime})$ is called \emph {strong isomorphism reducibility} of one pair to another. We will write $\mathsf{D}\leq\is\mathsf{D^\prime}$.
  If  $\mathsf{D}\leq\is \mathsf{D^\prime}$ and  $\mathsf{D^\prime}\leq\is \mathsf{D}$, i.e.,
   $\mathsf{D}\equiv\is \mathsf{D^\prime}$, then $\mathsf{D}$ and  $\mathsf{D^\prime}$  {\it have the same strong isomorphism degree} (see \cite {BChen}). It  was proven \cite {BChen} that the classes of finite sets, finite fields, finite abelian groups, finite cyclic groups, and finite sets with linear orderings all have the same strong isomorphism degree. However, as was also shown in \cite{BChen}, the problem of classifying undirected graphs is harder than the problem of  classifying all finite groups.

  In Section 5, we prove that the classification problem for the class of graphs is harder than the classification problem for  the class of graph Lie algebras and a class of finite $p$-groups.

 We also investigate a relation of the above classification problems to the well-known problem of classifying pairs of matrices over a field up to similarity. To be precise, let us denote by $\W_{1}$ the set of all pairs of $n\times n$ matrices, for all $n\in \mathbb{N}$, over a field $K$, and by $\W_{2}$ the set of all transformations of simultaneous similarity of pairs of matrices from $\W_{1}$: $$(A,B)\mapsto(S^{-1}AS, S^{-1}BS),$$
 where $A,B\in M(n,K)$ and $S\in \rm{GL}(n,K)$. This defines the pair $\W=(\W_{1},\W_{2})$. The classification problem for $\W$ ($\W$-problem)  is {\it the problem of classifying pairs of matrices up to similarity}. A matrix problem  is called  \emph{wild} if it contains the $\W$-problem as a subproblem. Wild problems are hopeless in a certain sense (see \cite{BS}).

Transformations from $\W_{2}$ induce an equivalence relation $T_{W}$  on $\W_{1}$. We say that the pair $(\W_{1}, T_{W})$
   corresponds to the pair $\W=(\W_{1},\W_{2})$. Let $\mathsf{\Omega}$ be a class of finite structures and $\rm{Iso}_\mathsf{\Omega}$
be the isomorphism relation on $\mathsf{\Omega}$. Let us fix a countable field $K$. The isomorphism problem for $\Omega$ is called \emph{Borel-wild} ($\mathcal{B}$-wild) over $K$ if the pair $(\W_1,T_{W})$ is polynomial time Borel-reducible to the pair  $(\Omega,\rm{Iso}_\Omega)$, i.e., $(\W_1,T_{W})\leq_B^P(\Omega,\rm{Iso}_\Omega)$. If $(\W_1,T_{W})\leq_B^P(\mathsf{\Omega}, \rm{Iso}_\mathsf{\Omega})$  but
$(\W_1,T_{W})\not\equiv_{\is}(\mathsf{\Omega},\rm{Iso}_\mathsf{\Omega})$, we say that the isomorphism problem for $\mathsf{\Omega}$ is \emph{Borel-superwild} and write $(\W_1,T_{W})<_P^B(\mathsf{\Omega}, \rm{Iso}_\mathsf{\Omega})$.

We prove that the class of undirected graphs without loops is Borel-superwild.
We also show that wildness of matrix problems over countable fields implies their Borel-wildness. The converse is an open problem.

  Below, all graphs are assumed to be finite undirected graphs without loops and multiple edges.

\section{\label{groups-section}A construction of a Lie algebras by a graph}
We give a reduction from the graph isomorphism
problem to the isomorphism problem for some class of $2$-nilpotent Lie algebras.

For each vector space $V$ over a field $K$, and a subset $W\subset V$, we denote by $\Span_K W$ the vector subspace of $V$ generated by all elements of $W$.

Denote by $F_n$ the free Lie algebra over $K$ generated by $u_{1},\dots,u_{n}$ and write $F_n^{3}:=\Span_K\{[[u_i,u_j],u_k]\,|\, i,j,k=1,\dots,n\}$. Then
\begin{equation}\label{freeni}
 N_n:=F_n/F_n^{3}
 \end{equation}
 is the free 2-step nilpotent algebra freely generated by $u_1+F_n^{3},\dots,u_n+F_n^{3}$.

Another realization of this algebra is given by M. Gauger \cite{G}. Let $V$ be the vector space over $K$ freely generated by $v_1,\dots,v_n$, and \[
\wedge^{2}V=V\wedge V:=V\otimes V/\Span_K\{v\otimes v\,|\, v\in V\}\] be the exterior square of $V$ (see \cite{Lang}). Turn the vector space $V\oplus \wedge^{2}V$ into a 2-step nilpotent Lie algebra in which the multiplication is given by
  \begin{equation}\label{metabelian}
[v_{i},v_{j}]=v_{i}\wedge v_{j},\quad [v_{i},v_{j}\wedge v_{k}]=[v_{i}\wedge v_{j},v_{k}]=[v_{i}\wedge v_{k},v_{j}\wedge v_l]=0,
\end{equation}
where $i,j,k,l=1,\dots,n.$ We identify $N$ and $V\oplus \wedge^{2}V$ via the isomorphism $\varphi: N\to V\oplus \wedge^{2}V$ that maps each $v_i\in N$ to $v_{i}\in V\oplus \wedge^{2}V$.

\begin{definition} \label{kdt}
Let $\Gamma=(T,E)$ be a graph with the vertex set $T=\{v_{1},\dots ,v_{n}\}$ and the edge set $E$.
Then for graph $\Gamma$ the following holds.
\begin{itemize}
  \item The {\it subspace of $N_n$ that corresponds to $\Gamma$} is the vector space
\[
I:=\Span_K\{v_{i}\wedge v_{j}\,|\,\{v_{i},v_{j}\} \in E\}\subseteq V\wedge V\subset N_n
\] (because the algebra $N_n$ defined in \eqref{freeni} is 2-step nilpotent, $I$ is an ideal of $N_n$).

  \item  The \emph{graph Lie algebra corresponding to $\Gamma$} is
 \[L(\Gamma):=N_n/I.\]
  \end{itemize}
  \end{definition}

In this section, we prove the following theorem.
\begin{theorem}\label{graphs-to-Lie algebra}
For any graphs  $\Gamma_{1}$ and $\Gamma_{2}$,
$$L(\Gamma_{1})\cong L(\Gamma_{2})\quad\Longleftrightarrow\quad
\Gamma_{1} \cong \Gamma_{2}.$$
\end{theorem}

 \begin{proof}
We use the following statements:
\begin{itemize}
  \item[(a)] Let $X=\{x_{1},\dots ,x_{n}\}$ and $Y=\{y_{1},\dots ,y_{n}\}$ be two bases of a vector space $V$ over a field $K$, where $\Gamma_{1}=(X,E_{1})$
  and $\Gamma_{2}=(Y,E_{2})$ are two graphs. Let $I_{1}$ and  $I_{2}$ be two subspaces of  $\wedge^{2}V$ corresponding to the graphs $\Gamma_{1}$ and $\Gamma_{2}$, respectively. The graphs $\Gamma_{1}$ and $\Gamma_{2}$ are isomorphic if and only if $I_{1}=I_{2}$ (see \cite{PBel}).

  \item[(b)]
 Let $N$ be a free  nilpotent Lie algebra of rank $n$. Then $N$ is freely generated by every system of $n$ generators of $N$ that are linearly independent modulo $N^{2}$ (see \cite{Mal_2}).

  \item[(c)]
  Let $L$ be a nilpotent Lie algebra and $\dim L/L^{2}=m$.  A subset $S=\{s_{1},\dots ,s_{m}\} \subseteq L$ generates $L$ if and only if the set $\{s+L^{2}|s\in S\}$ is a basis of $L/L^{2}$ (see \cite{Mal_2}).
  \end{itemize}

 Observe that if $\varphi:\Gamma_{1}\rightarrow\Gamma_{2}$ is any graph isomorphism, it induces the natural isomorphism
between $L_{1}=L(\Gamma_{1})$ and $L_{2}=L(\Gamma_{2})$. So we only have to prove that if
 \begin{equation}\label{nrf}
 \tau: L_{1}\rightarrow L_{2}
\end{equation}
is an isomorphism from $ L_{1}$ to  $L_{2}$,  $\Gamma_{1}\cong\Gamma_{2}$.

Let $N_{1}$ and $N_{2}$  be two free 2-step nilpotent Lie algebras generated by sets of vertices $T_{1}=\{v_{1},\dots ,v_{n_{1}}\}$
and $T_{2}=\{u_{1},\dots ,u_{n_{2}}\}$  of our graphs $\Gamma=(T_{1},E_{1})$ and $\Gamma=(T_{2},E_{2})$, respectively. Let $V$ and $U$ be the vector spaces over $K$ freely generated by the sets $T_{1}$ and $T_{2}$, respectively. Write
\begin{align*}
N_{1}&:=V\oplus \wedge^{2}V,& I_{1}&:=\Span _{K}\{v_{i}\wedge v_{j}| \{v_{i},v_{j}\} \in E_{1}\},\\
N_{2}&:=U\oplus\wedge^{2}U,&I_{2}&:=\Span _{K}\{u_{k}\wedge u_{m}| \{u_{k},u_{m}\} \in E_{2}\}.
\end{align*}
By the definition of graph Lie algebra we can write $ L_{1}=N_{1}/I_{1}$ and $ L_{2}=N_{2}/I_{2}$, where $I_1$ and $I_2$ are the vector spaces corresponding to the graphs $\Gamma_1$ and $\Gamma_2$. Because the algebras $L_{1}$ and $L_{2}$ are isomorphic, $L_{1}/L_{1}^{2}\cong L_{2}/L_{2}^{2}$. Using \[\mbox{dim}L_{i}/L_{i}^{2}=\dim N_{i}/N_{i}^{2}=n_{i},\qquad i=1,2,\]  we get  $n_1= n_2$. Write $n:=n_1=n_2$.

Consider the diagram
\begin{equation}\label{nh5}
\begin{split}
\xymatrix{N_{1}\ar[d]^{\pi_{1}}\ar@{-->}[r]^{\varphi} & N_{2}\ar[d]^{\pi_{2}}  \\
N_{1}/I_{1}\ar[r] ^{\tau} &N_{2}/I_{2}}
\end{split}
\end{equation}
where $\pi_{1}$ and $\pi_{2}$ are the canonical surjections and $\tau$ is the isomorphism \eqref{nrf}. Because $\tau\pi_{1}$ is a surjective map, there exist elements $w_1,\dots,w_n\in N_{2}$ such that $\pi_{2}(w_{i})=\tau\pi_{1}(v_{i})$ for all $i=1,\dots,n$.

Let us show that the elements $w_{1},\dots ,w_{n}$ are independent modulo $N_{2}^{2}$. For otherwise, elements $\tau^{-1}\pi_{2}(w_{1}),\dots,\tau^{-1}\pi_{2}(w_{n})$ are dependent modulo $L_{1}^{2}$. Therefore,
the elements $\pi_{1}(v_{1}),\dots,\pi_{1}(v_{n})$ are dependent modulo $L_{1}^{2}$.  As a consequence, the elements $v_1,\dots,v_n$ are dependent modulo $N_{1}^{2}$, which is impossible.

Let us define the homomorphism $\varphi:N_{1}\rightarrow N_{2}$ such that $\varphi(v_{i})= w_{i}$ for $i=1,\dots ,n$. By the statement (b),
the elements  $w_{1},\dots ,w_{n}$ freely generate the algebra $N_{2}$. Hence the homomorphism $\varphi$ is an isomorphism.
 Since $\tau\pi_{1}(v_{i})=\pi_{2}\varphi(v_{i})$ and the elements $v_1,\dots,v_n$  generate $N_{1}$,  diagram \eqref{nh5} is commutative. Therefore, $\varphi(I_{1})=I_{2}$.

 Let us write each $w_i$ as the sum
\[
w_{i}=\alpha_{i1}u_{1}+\dots+\alpha_{in}u_{n}+ b_{i},
\]
in which $b_{i}\in N_{2}^{2},\ u_{k}\in T_{2},\ \alpha_{ik}\in K,$ and                        $i,k=1,\dots ,n$. The elements $\varphi (v_{i}),i=1,\dots ,n$ generate the algebra $N_{2}$. By  statement (c), the elements \[ d_{i}:=\alpha_{i1}u_{1}+\dots+\alpha_{in}u_{n},\quad i=1,\dots ,n,\] are linearly independent modulo $N_{2}^{2}$ and they generate the algebra $N_{2}$.  By the statement (b), the elements $d_1,\dots,d_n$ freely generate $N_{2}$.
  Consider the homomorphism $\psi:N_{1}\rightarrow N_{2}$ such that
\begin{equation}\label{hye}
   \psi(v_{i})=d_{i},\qquad i=1,\dots ,n.
\end{equation}
The elements $d_{1},\dots ,d_{n}$ freely generate $N_{2}$, hence the homomorphism $\psi$ is an isomorphism.

  Because $N_{2}$ is a 2-step nilpotent Lie algebra and $\psi,\varphi:N_{1}\rightarrow N_{2}$ are Lie homomorphisms, we have
  $$ \varphi(v_{i}\wedge v_{j})=(d_{i}+b_{i})\wedge (d_{j}+b_{j})=d_{i}\wedge d_{j}=\psi(v_{i}\wedge v_{j}).$$
Therefore, $\psi(I_{1})=I_{2}$ and so $\psi(V)=U$ by \eqref{hye}.

Consider the graph $\Gamma_{3}=(T_{3},E_{3})$ with
\[
T_{3}:=\{d_1,\dots,d_n\},\quad E_3:=\{\{d_{i}, d_{j}\}\,|\,\{v_{i},v_{j}\}\in E_{1}\}.
\]
Denote by  $D$ the vector space freely generated by $T_{3}$ and consider the Lie algebra $N_{3}:=D\oplus \wedge^{2}D$. Denote by $I_3$ the subspace of $N_3$ that corresponds to $\Gamma_3$ (in the sense of Definition \ref{kdt}). The equalities
$d_{i}\wedge d_{j}=\psi(v_{i}\wedge v_{j})$ ensure $\psi(I_{1})=I_{3}$. Because  $\psi(I_{1})=I_{2}$, we get  $I_{3}=I_{2}$.
By the statement (a), the graphs $\Gamma_{3}$ and $\Gamma_{2}$ are isomorphic. Hence, the graphs $\Gamma_{1}$ and $\Gamma_{2}$ are also isomorphic.
    \end{proof}

\section{A construction of a $p$-group by a graph}\label{ksy}
 First we give a brief summary of the Lazard correspondence (\cite{Laz}, see also \cite{Eick}) between the category of nilpotent Lie rings $\mathcal{L}$ of nilpotency class $c$ and order $p^n$, $p>c$, and  category $\mathcal {G}$ of finite $p$-group of order $p^n$ and nilpotency class $c$. For each  $L\in \mathcal{L}$ we denote by $\Gr(L)\in \mathcal{G}$ the group with the same set of elements and with multiplication defined by the Beiker--Campbell--Hausdorff formula (BCH-formula,
\cite{Laz}), which has the form:
\begin{equation}
 g\cdot h=g+h+h_{1}(g,h)= g+h+\frac{1}{2}[g,h]+\frac{1}{12}[g,g,h]+\cdots,\quad g,h\in L,   \label{Lie_group}
\end{equation}
 where $h_{1}(g,h)$ is a finite linear combination over the field of rational numbers $\mathbb{Q}$ of Lie ring commutators in $g$ and $h$. The coefficients of the above linear combination are given as rationals whose denominators are not divisible by any prime greater than $c$. Thus the element $h_{1}(g,h)$ can be evaluated in $L$. Note that the expression for the element $h_{1}(g,h)$ only depends on the nilpotency class $c$, but not on $p$ or $L$ (for more details see \cite{Laz}).

 Conversely, let $G$ be a group in $\mathcal{G}$. Turn $G$ into the Lie algebra, in which the Lie operations $+$ and  $[\ ,\,]_{L}$ are defined as follows:
\begin{equation}  \label{Lie_algebra}
 g+h:=g\cdot h\cdot h_{2}(g,h),\qquad
 [g,h]_{L}:=[g,h]_{G} \cdot h_{3}(g,h).
\end{equation}
Here $g,h\in G$; $[g,h]_{G}=g^{-1}h^{-1}gh$ is the group commutator; and $h_{2}(g,h)$ and $h_{3}(g,h)$ are (defined in \cite{Laz}) products of formal powers of the group commutators of $g$ and $h$ (the expressions \eqref{Lie_algebra} are called the \emph{inverse  $\rm{BSH}$-formulas}). Because the denominators of exponents in the expressions of $h_{2}(g,h)$ and $h_{3}(g,h)$ are not divisible by any prime greater than $c$, they can be evaluated in finite $p$-group $G$. Denote the above Lie ring by $\Lie(G)$. Note that the expressions for the elements $h_{1}(g,h)$ and  $h_{2}(g,h)$ only depend on the nilpotency class $c$, but not on $p$ or $G$.

It can be proved that  $\Gr(\Lie(G))=G$ and $\Lie(\Gr(L))=L$ hold for a group $G\in\mathcal{G}$ and a Lie ring $L\in\mathcal{L}$. We say that $G$ and $L$ are Lazard correspondening to each other. The Lazard correspondence also gives an isomorphism between  category $\mathcal{L}$ of nilpotent Lie rings of order $p^{n}$ and nilpotency class $c$, and category $\mathcal{G}$ of finite $p$-groups of nilpotency class $c$, provided $p<c$.

In this section, we use the Lazard correspondence to describe a relation between the isomorphism problem for graphs and for a class of $p$-groups corresponding to them. Let $L$ be a finite dimensional Lie algebra over the field $\mathbb{F}_{p}=\mathbb{Z}/p\mathbb{Z}$ with $p\neq 2$. In what follows, we denote by $L^R$ the Lie ring of a Lie algebra $L$. It is evident that two finite dimensional Lie algebras, $L_{1}$ and $L_{2}$, over the field $\mathbb{F}_{p}$ are isomorphic if and only if $L^{R}_{1}$ is isomorphic to  $L^{R}_{2}$.

 Let $\Gamma=(T,E)$ be a graph. As in Section \ref{graphs-to-Lie algebra}, we define the vector space $V$ freely generated by the set of vertices $T=\{v_1,\dots,v_n\}$ over finite field $\mathbb{F}_{p}$ with $p\neq 2$ and the free 2-step nilpotent Lie algebra $N_n=V\bigoplus \wedge^{2}V$ with defining relation (\ref{metabelian}). Let $M= \Gr(N_n^{R})$ be the group Lazard corresponding to  Lie ring $N_n^{R}$. Because $N_n^{R}$ is a 2-step nilpotent Lie ring of characteristic $p\neq 2$, a multiplication on $M$ can be defined by BCH-formula (\ref{Lie_group}):
\begin{equation} \label{met_2}
(v_{1}+w_{1})(v_{2}+w_{2})= v_{1}+v_{2}+w_{1}+w_{2}+1/2 (v_{1}\wedge v_{2}),
\end{equation}
for all $ v_{1}, v_{2}\in V$ and $ w_{1}, w_{2}\in \wedge^{2}V$. Note that $N^{R}$ is a free ring in the variety of Lie rings determined by the identities: $ p\cdot x=0, [[x,y],z]=0$. Because any Lie ring homomorphism $\varphi:L^{R}_{1}\rightarrow L^{R}_{2}$, where $L^{R}_{1},L^{R}_{2}\in\mathcal{L}$, induces the group homomorphism
 $\hat\varphi:\Gr(L^{R}_{1})\rightarrow \Gr(L^{R}_{2})$, $M_n$ is a free group freely generated by $T$ in the variety of groups determined by the identities $ x^{p}=1, [[x,y],z]=1$ (see formula (\ref{Lie_group})).

 As in the case of Lie algebras, we can define a 2-step nilpotent finite $p$-group corresponding to the graph $\Gamma=(T,E)$.
  \begin{definition}
  For each graph $\Gamma=(T,E)$, define
  \begin{itemize}
    \item the subgroup of $J$ of $M_n$ generated by $v_{i}\wedge v_{j}$, where $\{v_{i},v_{j}\} \in E$ (because the group $M_n$ is  2-step nilpotent, $J$ is a normal subgroup of $M_n$);
    \item \emph{the graph $p$-group corresponding} to $\Gamma$ is
     $$G(\Gamma)=M_n/J,$$

      \end{itemize}
      \end{definition}

  Below we need the following result:
  \begin{proposition} [{\cite{Eick}},{\cite{Laz}}]\label{Laz}
 {\it Let $G$ be a finite $p$-group of class $c<p$, and  $H$ be its Lazard correspondent ring. Let $G_{0}$ be a normal subgroup in $G$ and $H_{0}$ be the corresponding ideal in $H$. Then
  $\psi:G/G_{0}\rightarrow H/H_{0} : xG_{0}\rightarrow x+H_{0}$ is a  well-defined bijection, and it induces the Lazard correspondence between $G/G_{0}$ and $ H/H_{0} $.} \qed
  \end{proposition}
\begin{proposition}\label{Laz_1}
Let $\Gamma=(T,E)$ be a graph, $G(\Gamma)$ be the graph $p$-group and $L(\Gamma)$ be the Lie graph algebra over the field $\mathbb{F}_{p}$ corresponding to the graph $\Gamma$. Then $G(\Gamma)$ and $L(\Gamma)$ are Lazard correspondents of each other.
 \end{proposition}
 \begin{proof}
 Let $I$ be the subspace of  ring $N_n$ corresponding to $\Gamma$. Then $I^R$ is an ideal of $N_n^R$ and $L(\Gamma)^R=N_n^R/I^R$. It is clear that  ring $N_n^R$ and  group $M_n$ are in Lazard correspondence. The normal subgroup $J$ and  ideal $I^R$ of  ring $N_n$ are also in Lazard correspondence, i.e., $\Gr(I^R)=J$. By Proposition \ref{Laz}  Lie ring $L(\Gamma)^{R}=N_n^{R}/I^R$ is the Lazard correspondent of group  $G(\Gamma)=M_n/J$, i.e., $\Gr(L(\Gamma)^R)=G(\Gamma)$.
 \end{proof}
  \begin{theorem}\label{graph-to-p-group}
For every two graphs  $\Gamma_{1}$ and $\Gamma_{2}$
$$G(\Gamma_{1})\cong G(\Gamma_{2})\quad\Longleftrightarrow\quad
\Gamma_{1} \cong \Gamma_{2}.$$
\end{theorem}
\begin{proof}
 If $G(\Gamma_{1})$ and $G(\Gamma_{2})$ are two isomorphic graph $p$-groups, their Lazard correspondent graph Lie rings $L(\Gamma_{1})^{R}$ and $L(\Gamma_{2})^{R}$ are also isomorphic. Hence  graph Lie algebras $L(\Gamma_{1})$ and $L(\Gamma_{2})$ over  field $\mathbb{F}_{p}$ are also isomorphic.  By Theorem \ref{graphs-to-Lie algebra}  graphs $\Gamma_{1}$ and $\Gamma_{2}$ are isomorphic.
The converse is trivial.
\end{proof}

\begin{remark}
 Theorem \ref{graph-to-p-group} can be proved  also using the known properties of locally finite varieties of $p$-groups (see \cite{PBel}). However, our proof reveals an important relation between graph Lie algebras and graph $p$-groups via the Lazard correspondence.
 \end{remark}

\section{Borel reducibility and wildness}
 We use Borel reducibility to define Borel-wildness ($\mathcal{B}$-wildness) of the isomorphism problem for classes of finite structures. Let $A$ be a countable set. Denote by $\Sigma$ a finite alphabet and by $\Sigma^{*}$ the free monoid over the alphabet $\Sigma$. As usual, a {\it language} over $\Sigma$ is a subset of the monoid  $\Sigma^{*}$. Encoding elements of $A$ by words from $\Sigma^{*}$  (this encoding can be done in many reasonable ways), we define a language $L_{A}$ over $\Sigma$.

 Let $R$ be an equivalence relation on $A$. The relation $R$ can be encoded as a language by taking the pairwise encoding of each pair in $R$. Hereinafter we  will abuse notation and write $(a,c)\in R$ (or $aRc$), where $a,c\in L_{A}$, for the equivalence relation $R$ on A, but what we really mean is $(a,c)\in L_{R}$,  where $L_{R}$ is the language over the alphabet $\Sigma$ induced by $R$.

 Let $A$ and $C$ be two countable sets. In the following we say that a map
 $f:A\rightarrow C$ is computable if the induced map $\hat f:L_{A}\rightarrow  L_{C}$ is computable.
  \begin{definition} \rm{\bf{\cite{Cos}}}
{\it Let $A$ and $C$ be two countable sets and $R,S$ be equivalence relations on $A$ and $C$, respectively. We say that $(A,R)$ is computably Borel-reducible to $(C,S)$, and write $(A,R)\leq_{B}(C,S)$, if there exists a computable map $f:A\rightarrow C$ such that for all
  $x$ and $y$ in $A$}
    $$
    xRy\Leftrightarrow f(x) S f(y).
    $$
 \end{definition}
  In other words, the reduction function $f$ yields a classification of the elements  of $A$ up to $R$ using invariants from $C/S$. We also say that $(A,R)$  and $(C,S)$ are \emph{Borel equivalent} and write $(A, R)\equiv_B (C, S)$ if they are Borel-reducible one to another, i.e., $(A,T)\leq_{B}(C,S)$ and $(C,S)\leq_{B}(A,R)$. If $f$ is computable in polynomial time, we say that $(A,R)$ is {\it polynomial-time Borel-reducible} to $(C,S)$ ($P$-Borel-reducible) and use the notation $(A,R)\leq_{B}^{P}(C,S)$. Similarly,
 we define $(A, R)\equiv^{P}_B (C,S)$.

 Let $\A_{1}$ be a set of $a$-tuples of matrices over a field $K$ and $\A_{2}$ be a set of admissible matrix transformations with them. Denote $\A=(\A_{1},\A_{2})$. The transformations from $\A_{2}$ induce the equivalence relation $T_\A$ on  set $\A_{1}$. The {\it classification problem} for the pair $\A=(\A_{1},\A_{2})$ is to find a description of the set of canonical $a$-tuples in the equivalence classes of the quotient set $\A_{1}/T_{\A}$. Hereafter, the classification problem for the pair $\A=(\A_{1},\A_{2})$ is called {\it an $\A$-matrix problem} (or shortly, an  {\it $\A$-problem}), (see \cite {BS}).
\begin{definition}[\cite{BS}]\label{matprob}
Given two pairs $\A=(\A_{1},\A_{2})$ and $\B=(\B_{1},\B_{2})$, we say that the $\A$-problem
\emph{is contained in} the $\B$-problem, ($\A\preceq \B $), if there exists a $b$-tuple $\T(x)=\T(x_{1},\dots ,x_{a})$ of matrices, whose entries are
non-commutative polynomials in $x_{1},\dots ,x_{a}$, such that
\begin{itemize}
\item $\T(A)=\T(A_{1},\dots ,A_{a})\in \B_{1}$ if $A=(A_{1},\dots ,A_{a})\in \A_{1}$;
\item for every $A,A\p\in \A_{1}$, $A$ reduces to $A\p$ by transformations from
$\A_{2}$ if and only if $\T(A)$ reduces to $\T(A\p)$ by
transformations from $\B_{2}$.
\end{itemize}
\end{definition}

If $\A\preceq \B $ and $\B\preceq \A $ we say that $\A= \B $. In this case a solution of the classification problem for $\B$ implies a solution of the clasification problem for $\A$.

Let us consider the pair $\T=(\T_{1}, T_{2})$, where $\T_{1}$ is the set of all square matrices of the order $n\times n$, for all $n\in \mathbb{N}$, over a field $K$, and where $\T_{2}$ is the set of all transformations of similarity of matrices from $\T_{1}$:
$$
 A\rightarrow S^{-1} A S,
$$
where $ A\in M(n,K)$ and $S\in \rm{GL}(n, K)$. A solution of the classification problem for  $\T$ over an algebraically closed field $K$ is  the canonical Jordan form of matrices from $\T_{1}$ (see \cite{Lang}).

Recall that the $\W$-problem defined in the Introduction is the classification problem for the pair $(\W_1, \W_2)$,
where $\W_1$ is the set of all pairs of $n\times n$ matrices, for all $n\in \mathbb{N}$, over a field $K$, and  $\W_{2}$ is the set of all transformations of simultaneous similarity of pairs of matrices from $\W_{1}$.  It can be proved that the $\W$-problem over an algebraically closed field $P$ strictly contains the $\T$-problem, i.e, $\T\preceq \W$ and $\T\neq \W$.

 In the notation of  definition \ref{matprob}, the admissible matrix transformations $\A_{2}$ (resp. $\B_{2}$) on $\A_{1}$ (resp. $\B_{1}$) define the equivalence relations $T_{\A}$ on $\A_{1}$ (resp. $T_{\B}$ on $\B_{1}$).
 \begin{proposition}\label{g-wild}
Let $K$ be a countable field and $\A=(\A_{1},\A_{2})$ and $\B=(\B_{1},\B_{2})$ be two pairs over $K$.  If the $\A$-problem is contained in the $\B$-problem, i.e., $\A\preceq \B$, then $(\A_{1}, T_{\A})\leq_{B}^{P}(\B_{1},T_{\B})$.

 \end{proposition}
 \begin{proof}
Let us fix an alphabet $\Sigma$ that contains all the symbols necessary to encode elements of  field $K$ and two additional symbols $\mid$ and $\parallel$. A matrix tuple is represented by words
from $\Sigma^{*}$, where the rows of a matrix are separated by $\mid$ and different matrices of the tuple by $\parallel$. This defines the languages $L_{\A_{1}}$ and $L_{\B_{1}}$ over the alphabet $\Sigma$. Using the $b$-tuple $\T(x)$ of matrices whose entries are non-commutative polynomials in $x_{1},\dots ,x_{a}$ we can construct mapping from $L_{\A_{1}}$ to $L_{\B_{1}}$ that is computable in  polynomial time. Hence, the pair $(\A_{1}, T_{\A})$ is $P$-Borel-reducible to $(\B_{1},T_{\B})$.

 \end{proof}
  \begin{definition}\rm{\bf{\cite{BS}}}\label{wild_1}
The classification  problem for  the pair $\A=(\A_{1},\A_{2})$ is called  wild if  the $\A$-problem contains the $\W$-problem, i.e., $\W\preceq \A$.
\end{definition}
The classification  problem for $\W$ is considered as hopeless in a certain sense. A list of some known wild matrix problems is given in \cite{BS}.

We now present another approach to the notion of wildness of a matrix problem over a countable field $K$. Let $A$ be a countable set and $R$ be an equivalence relation on $A$. The classification problem for the pair $(A,R)$ is to find a description of the set of canonical representatives in the equivalence classes of the quotient set $A/R$. To characterize the complexity of this classification problem we compare it to the complexity of the $\W$-problem over field $K$.

\begin{definition}\label{wild_3}
The classification  problem for  the pair $(A,R)$ is called Borel-wild ($\mathcal{B}$-wild) over $K$ if the pair  $(\W_1,T_{W})$ is $P$-Borel-reducible to the pair $(A,R)$, i.e., $(\W_1,T_{W})\leq_{B}^{P}(A,R)$.
\end{definition}
Let  $\A$-matrix problem be determined by the pair $(\A_{1},\A_{2})$ over a countable field $K$ and $\A^{\prime}=(\A_{1},T_{\A})$ be the pair corresponding to $(\A_{1},\A_{2})$. From Proposition \ref{g-wild} follows\ that if the $\A$-problem is wild, the pair  $\A^{\prime}=(\A_{1},T_{\A})$ is $\mathcal{B}$-wild. An interesting open question is whether the converse is also true.

We now define a notion of  Borel-wildness of the isomorphism problem for classes of finite structures. Let us recall the definition of \emph{a structure} (see \cite{Ebb, Mal_1}). A signature (or vocabulary) $\sigma$  is a finite sequence of relation symbols, function symbols, and constant symbols. Then, a structure $\rm{S}$ over the signature  $\sigma$ is defined as a tuple that includes an universe $U_{\rm{S}}$  and  an interpretation of all symbols from $\sigma$, i.e., an assignment of meaning to the symbols from $\sigma$ in $U_{\rm{S}}$. A structure $\rm{S}$ is finite if its universe $U_{\rm{S}}$ is finite.
The cardinality of the universe $U_{\rm{S}}$ will be denoted by $|U_{\rm{S}}|$.

{\it From this point forward we will work only with classes of  finite structures.}

Let $\mathsf{D}$ be a class of structures. Let us fix a finite alphabet $\Sigma$. We now encode a structure $T$ in $\mathsf{D}$ by words from $\Sigma^{*}$, $\rm{enc}(T)$. Denote  $L_\mathsf{D}=\{enc(A)|A\in\mathsf{D}\}$. We assume that the mappings $A\mapsto \rm{enc}(A)$ and $\rm{enc}(A)\mapsto A$ are computable in polynomial time. Let  $\mathsf{D}$ and $\mathsf{D^\prime}$ be two structures. We say that a map $f:\mathsf{D}\rightarrow \mathsf{D^\prime}$ is computable if the induced map $\hat f:L_\mathsf{D}\rightarrow  L_\mathsf{D^\prime}$ is computable.
\begin{definition}[{\cite{BChen}}, see also {\cite{Cos}}]
 Let $\mathsf{D}$ (resp. $\mathsf{D^\prime}$) be two classes of structures. We say that $\mathsf{D}$ is \emph{ strongly isomorphism-reducible} to $\mathsf{D^\prime}$, and write
  $\mathsf{D}\leq\is\mathsf{D^\prime}$, if there exists a function $f:\mathsf{D}\rightarrow \mathsf{D^\prime}$ computable in polynomial time and such that for all $A$ and $B$ in $\mathsf{D}$
 $$
 A\cong B\Leftrightarrow f(A)\cong f(B)
 $$
 \end{definition}
 If  $\mathsf{D}\leq\is\mathsf{D^\prime}$ and  $\mathsf{D^\prime}\leq\is \mathsf{D}$, $\mathsf{D}$ and $\mathsf{D^\prime}$  {\it have the same strong isomorphism degree}; we write $\mathsf{D}\equiv\is \mathsf{D^\prime}$. The equivalence $\equiv\is$ we will call  {\it $\rm{SID}$-equivalence}.

 Denote by $\rm{Iso}_\mathsf{D}$ (resp. $\rm{Iso}_\mathsf{D^\prime}$) the isomorphism relations on classes $\mathsf{D}$ and $\mathsf{D^\prime}$, respectively. It is clear that  $\mathsf{D}$ is  strongly isomorphism-reducible to $\mathsf{D^\prime}$ if and only if the pair ($\mathsf{D}, \rm{Iso}_{\textsc D}$) is $P$-Borel-reducible to ($\mathsf{D^\prime}, \rm{Iso}_\mathsf{D^\prime}$).

Let $K$ be a finite field and  $\W=(\W_{1},\W_{2})$ be the aforementioned pair over $K$.
\begin{definition}\label{wild_2}
The isomorphism problem for $\Omega$ is called Borel-wild ($\mathcal{B}$-wild) over $K$ if the pair $(\W_1,T_{W})$ is $P$-Borel-reducible to the pair  $(\Omega, \rm{Iso}_\Omega)$, i.e., $(\W_1,T_{W}) \leq_{B}^{P}(\Omega, \rm{Iso}_\Omega)$.
 \end{definition}

\begin{definition}\label{superwild}
 We say that the isomorphism problem for $\Omega$ is \emph{Borel-superwild} and write $(\W_1,T_{W}) <_{B}^{P}(\Omega, \rm{Iso}_\Omega)$, if it is Borel-wild and
 $(\W_1,T_{W})\not\equiv_{\is} (\Omega,\rm{Iso}_\Omega)$,
 \end{definition}
 In what follows, we omit the sign of the isomorphism relation defined on the $\Omega$ class and write  $(\W_1,T_{W})<_P^B\mathsf{\Omega}$.

  Let $P$ be a field of characteristic different from $2$. It is known that the isomorphism problems are wild for the following classes:
 \begin{itemize}
 \item finite dimensional Lie algebras over $P$ with cenral commutator subalgebra of dimension $3$ ( see \cite{GBel, BDLST}),
 \item local commutative associative algebras over $P$ with zero cube radical (see \cite{BBLPS}),
  \item finite $p$-groups of exponent $p$ with central commutator subgroup of order $p^{3}$ (see \cite{S}).
 \end{itemize}
 Note that wildness of the isomorphism problems for the first two classes means  wildness of the corresponding matrix problems in the sense of Definition \ref{wild_1}. However, wildness of the isomorphism problem for the third class of  finite $p$-groups should be understood in the sense of Definition \ref{wild_2}, where $P=\mathbb{F}_{p}$, i.e., as  Borel-wildness over $\mathbb{F}_{p}$ (see \cite{S}) . We will use these results to show Borel-wildness of the isomorphism problems for several classes of finite structures.

  \section{\label{comp} The complexity of the isomorphism problems}
  In the previous sections we have proved that the isomorphism problem for the  class of undirected graphs, denoted by  $\mathsf{GRAPH}$, can be reduced to the isomorphism problems for the class $\mathsf{GLA}$ of graph Lie algebras over the field $\mathbb{F}_{p}$ with $p\neq 2$, and the class $\mathsf{ GpG}$ of graph $p$-groups with $p\neq 2$, and vice versa. Now we prove that the isomorphism problem for $\mathsf{GRAPH}$ is harder than the isomorphism problems for the classes $\mathsf{GLA}$ and $\mathsf{GpG}$, and  is superwild. First we show that the isomorphism problems for the classes $\mathsf{ GpG}$ and $\mathsf{GLA}$ have the same isomorphism degree:
  \begin{theorem}\label{F}
        $\mathsf{GLA}\equiv\iso\mathsf{GpG}$
   \end{theorem}
   \begin{proof}
    The map $g:\mathsf{GLA}\rightarrow \mathsf{GpG}$ and  the map $f:\mathsf{GpG} \rightarrow \mathsf{GLA}$ realizes the Lazard correspondence between the classes  $\mathsf{GLA}$ and $\mathsf{GpG}$, are polynomially computable (see the formula (\ref{met_2}) and the inverse BCH-formulas (\ref{Lie_algebra})). Because maps $f$ and $g$ are isomorphism-preserving, we have $\mathsf{GLA}\equiv\is \mathsf{GpG}$.
\end{proof}

 The following result was proved in \cite{Ebb} (see also {\cite{BChen}}):
  \begin{proposition}\label{E}
 $\mathsf{\textsc A}\leq\iso \mathsf{GRAPH}$, for any class of structures {\textsc A}.
\end{proposition}

\begin{theorem}\label{s2}
 The isomorphism problem for the class $\mathsf{GRAPH}$ is superwild,
   $$(\W_1,T_{W})\prec_{B}^{P}\mathsf{GRAPH},$$
   and is harder than the isomorphism problems for the classes $\mathsf{GLA}$ and $\mathsf{GpG}$, i.e.,
 $$\mathsf{GpG}\prec_{\iso}\mathsf{GRAPH},\ \ \ \mathsf{GLA}\prec_{\iso}\mathsf{GRAPH}.$$
 \end{theorem}

\begin{proof}
In \cite {S} was proven that the isomorphism problem for the class of finite $p$-groups with a central commutator subgroup of order $p^{2}$ is wild.  By Proposition \ref{g-wild} it is $\B$-wild. Therefore, the isomorphism problem for the class of finite groups, denoted by $\mathsf{GROUP}$, is $\B$-wild, i.e., $(\W_1,T_{W})\preceq_{B}^{P}\mathsf{GROUP}$. It is known \cite{BChen} that $\mathsf{GROUP}\prec\is\mathsf{GRAPH}$. Therefore, $(\W_1,T_{W})\prec_{B}^{P}\mathsf{GRAPH}$, i.e., the isomorphism problem for the class of graph is superwild.

Because $\mathsf{Group}\prec\is\mathsf{Graph}$, we have $\mathsf{GpG}\prec\is\mathsf{GRAPH}$. By  Theorem
\ref{E}, $\mathsf{GLA}\equiv\is \mathsf{GpG}$. Hence, $\mathsf{GLA}\prec\is\mathsf{GRAPH}$.

The last relation can also be proven directly in the context of the theory of Lie algebras.
Indeed, according to Proposition \ref{E} we have $\mathsf{GLA}\preceq_{\is}\mathsf{GRAPH}$. Let us show that
$\mathsf{GLA}\not\equiv_{\is}\mathsf{GRAPH}$. Since $\mathsf{GLA}\preceq_{\is} \mathsf{GRAPH}$, there exists a computable function  $f:\mathsf{GLA}\rightarrow \mathsf{GRAPH}$.  Hence, there exists a polynomial $g(x)$ such that for a Lie algebra $L\in\mathsf{GLA}$,  $|f(L)|\leq g(|L|)$, where $|L|$ denotes the cardinality of the algebra $L$.  Because $f$ is a strong isomorphism reduction, the number of isomorphism types $N_{1}$ of nilpotent Lie algebras of the cardinality $\leq p^m$ over field $\mathbb{F}_{p}$ is equal to the number of isomorphism types $N_{2}$  of graphs with a number of vertices $\leq g( p^m)$. It is known
\cite {Black} that the number of nilpotent Lie algebras with $p^m$ elements over the field $\mathbb{F}_{p}$ is at most $p^{\frac{2}{27} m^3+O(m^{5/2})}$. Hence, $N_1\leq mp^{\varphi(m)}$, where $\varphi(x)$ is a polynomial. On the other hand $N_2\geq 2^{\frac{1}{2}g(p^m)(g(p^m)-1)}$. Therefore, for a sufficiently large $m$, $N_1< N_2$. We arrive at a contradiction. Hence, we obtain again that $\mathsf{GLA}\prec\is\mathsf{GRAPH}$.
\end{proof}
For the classes $\mathsf{GLA}$ and $\mathsf{GpG}$ we can show more than $\rm{SID}$-equivalence. We need the following definitions.
\begin{definition}[{\cite{Z}}]
  Let $\M$ and $\T$ be two categories with the classes of objects $\rm{Ob}(\M)$ and $\rm{Ob}(\T)$, respectively. Let  $\varphi:\M\rightarrow\T$ be a functor. If for any objects $X,Y$ from $\rm{Ob}(\M)$ the induced mapping $\varphi^\prime: \rm{Mor}(X,Y)\rightarrow \rm{Mor}(\varphi(X),\varphi(Y))$ is bijective, then $\varphi$ is called a complete embedding of $\M$ into $\T$.
  \end{definition}

 If there exists a complete embedding $\varphi$ of a category $\M$ into a category $\T$ such that the induced mapping $\varphi'':\rm{Ob}(\M)\rightarrow\rm{Ob}(\T)$ is of a polynomial complexity, then the isomorphism problem for the class of objects $\rm{Ob}(\M)$  functorially reduces to the analogous problem for the class $\rm{Ob}(\T)$. We say that the isomorphism problem for the class  $\rm{Ob}(\M)$ is functorialy equivalent to the same problem for the class $\rm{O}b(\T)$ if they are functorially reducible one to another.

  Let us regard two pairs $(\rm{Ob}(\M), Iso_{\M})$ and $(\rm{Ob}(\T), Iso_{\T})$, where $ \Iso_{\M}$ (resp. $\Iso_{\T}$) denotes the isomorphism relation on the class $\rm{Ob}(\M)$ (resp. $\rm{Ob}(\T)$). A functorial reduction of the isomorphism problem for the objects from $\rm{Ob}(\M)$  to the analogous problem for $\rm{Ob}(\T)$  is more restrictive than a  $P$-Borel reduction of the pair $(\rm{Ob}(\M), Iso_{\M})$
to $(\rm{Ob}(\T), Iso_{\T})$, because it requires existense of a bijection between the sets of isomorphisms of the two categories. Therefore, a functorial reduction implies a  $P$-Borel reduction of the above  pairs (or a strong isomorphism reduction of the class $\rm{Ob}(\M)$ to $\rm{Ob}(\T)$).

 Let us regard the above mentioned classes of finite structures as categories. For conveniences we designate these categories by the same letters as the corresponding classes. Then the isomorphism problems for the objects of the category $\mathsf{GLA}$ and the category $\mathsf{GpG}$ are functorially equivalent. This follows from the following property of Lazard correspondence of these categories: any Lie  isomorphism $\psi:L_{1}\rightarrow L_{2}$, where the Lie rings $L_{1}$ and $L_{2}$ belong to the class $\mathsf{GLA}$, induces a group isomorphism $\hat\psi:\Gr(L_{1})\rightarrow \Gr(L_{2})$ of the corresponding graph $p$-groups $\Gr(L_{1})$ and $ \Gr(L_{2})$ from the class $\mathsf{GpG}$, and vice versa.

   \section{Acknowledgments}
The authors are grateful to G. Belitskii and V. Sergeichuk for fruitful discussions and interest in this work.




\begin{thebibliography}{00}
\bibitem{PBel} P. Beletskii, A class of finite $p$-groups,
{\it Problems in group theory and homological algebra} {\bf 163} (1985) 3--13.

\bibitem{GBel} G. Belitskii, R. Lipyanski and V. Sergeichuk,  Problems of classifying associative or Lie algebras and triples
 of symmetric or skew-symmetric matrices are wild, {\it Linear Algebra Appl.} {\bf 407} (2005) 249--262.

 \bibitem{BDLST}  G. Belitskii, A. Dmytryshyn, R. Lipyanski, V. Sergeichuk and A. Tsurkov, Problems of classifying associative or Lie algebras  over a field of characteristic not $2$
 and finite metabelian groups are wild, {\it Electron. J. Linear Algebra} {\bf 18} (2009) 516--529.
\bibitem{BS} G. Belitskii and V. Sergeichuk, Complexity of matrix problems, {\it Linear Algebra Appl.} {\bf 361} (2003) 203--222.
\bibitem{BBLPS} G. Belitskii, V. Bondarenko, R. Lipyanski, V. Plachotnik and V. Sergeichuk,  The problems of classifying pairs of forms and local algebras with zero cube radical are wild,
  {\it Linear Algebra Appl.} {\bf 402} (2005) 135--142.
  \bibitem{Black} S. Blackburn, P. M. Neumann, G. Venkataraman, {\it Enumeration of finite groups} (Cambridge University Press, Cambridge, 2007)
  \bibitem{BChen}   S. Buss, Y. Chen, J. Flum, S. Fridman and M. Muller, Strong isomorphism reductions in complexity theory,{\it J. Symbolic Logic } {\bf 76} (4) (2011) 1381--1402.
  \bibitem{Cos} S. Coskey, J.D. Hamkins and R. Miller, The hierarchy of equivalence relations on the natural numbers under computable reducibility, {\it Computability } {\bf 1} (1) (2012) 15--38.
\bibitem{D} C. Droms, Isomorphism of graph groups, {\it Proc. Amer. Math. Soc.} {\bf 100} (3) (1987) 407--409.
\bibitem{Ebb}H.-D. Ebbinghaus and J. Flum, {\it Finite model theory} (Springer-Verlag, Berlin, 1999).
\bibitem{Eick} B. Eick, M. Horn and S. Zandi, Shur multipliers and the Lazard correcpondence, {\it Arch. Math.} {\bf 90} (2012) 217--226.
 \bibitem{F}H. Friedman and L. Stanley,  Borel reducibility theory for classes of countable structures, {\it J. Symbolic Logic} {\bf 54} (3) (1989) 894--914.
 \bibitem{G} M. Gauger, On the classification of metabelian Lie algebras, {\it Trans. Amer. Math. Soc.} {\bf 179} (1973) 293--329.
 \bibitem{Kim} K.H. Kim, L. Makar-Limanov, J. Neggers, F.W. Roush,  Graph algebras. {\it J. Algebra} {\bf 64} (1) (1980) 46--51.
 \bibitem {Laz} M. Lazard, Sur les groupes nilpotents et les anneaux de Lie,  {\it Ann. Sci. Ecole Norm. Sup.} {\bf 71} (3) (1954) 101--190.
 \bibitem{Lang} S. Lang, {\it Algebra} (Springer-Verlag, New York, 2002).
  \bibitem {Mal_1} A.I. Mal'tsev,  {\it Algebraic systems} (Springer-Verlag, New York-Heidelberg, 1973).
  \bibitem {Mal_2} A.I. Mal'tsev, On algebras defined by indentities, {\it Mat.Sb.} {\bf 26} 1950 19--23 (Ruissan).
  \bibitem {Mil} G. Miller, Graph isomorphism, general remarks, {\it J. Comput. System Sci. } {\bf 18} (2) (1979) 128--142
  \bibitem {Sax} N. Kayal, N. Saxena, Complexity of ring morphism problems,  {\it J. Comput. Complexity } {\bf 15} (4) (2006) 342--390
  \bibitem {S} V. Sergeichuk, The classification of metabelian $p$-groups, {\it in: Matrix problems, Akad. Nauk Ukrain. SSR Inst. Mat., Kiev}, (1977) 150--16 (Russian).
   \bibitem{Z}V. Zemlyachenko, N. Korneenko and R. Tyshkevich, The graph isomorphism problem,  {\it Zap. Nauchn. Sem. Leningrad. Otdel. Mat. Inst. Steklov} {\bf 118} (1982) 83--158.


\end{thebibliography}
\end{document}